# Non- and semi-parametric analysis of failure time data with missing failure indicators[*]

## Irene Gijbels[1],[†], Danyu Lin[2],[‡] and Zhiliang Ying[3],[§]

*Katholieke Universiteit Leuven, University of North Carolina and Columbia University*

**Abstract:** A class of estimating functions is introduced for the regression parameter of the Cox proportional hazards model to allow unknown failure statuses on some study subjects. The consistency and asymptotic normality of the resulting estimators are established under mild conditions. An adaptive estimator which achieves the minimum variance-covariance bound of the class is constructed. Numerical studies demonstrate that the asymptotic approximations are adequate for practical use and that the efficiency gain of the adaptive estimator over the complete-case analysis can be quite substantial. Similar methods are also developed for the nonparametric estimation of the survival function of a homogeneous population and for the estimation of the cumulative baseline hazard function under the Cox model.

## 1. Introduction

Let $(T_i, C_i, Z_i')$ $(i = 1, \ldots, n)$ be $n$ independent replicates of the random vector $(T, C, Z')$, where $T$ and $C$ denote the failure and censoring times, and $Z$ denotes a $p \times 1$ vector of possibly time-varying covariates. The observations consist of $(X_i, \delta_i, Z_i')$ $(i = 1, \ldots, n)$, where $X_i = T_i \wedge C_i$ and $\delta_i = 1_{(T_i \leq C_i)}$. Assume that $T_i$ and $C_i$ are conditionally independent given $Z_i$.

The widely-used Cox semiparametric regression model [4] postulates that, conditional on $Z(t)$, the hazard function $\lambda(t)$ for $T$ takes the form $e^{\beta_0' Z(t)} \lambda_0(t)$, where $\beta_0$ is a $p$-dimensional regression parameter and $\lambda_0(\cdot)$ is an unspecified baseline hazard function. The maximum partial likelihood estimator $\hat{\beta}_f$ for $\beta_0$ is obtained by maximizing

$$(1.1) \qquad L(\beta) = \prod_{i=1}^{n} \left\{ \frac{e^{\beta' Z_i(X_i)}}{\sum_{j=1}^{n} 1_{(X_j \geq X_i)} e^{\beta' Z_j(X_i)}} \right\}^{\delta_i},$$

---

[*]Part of the work was done while the first and the third authors were visiting the Mathematical Sciences Research Institute, Berkeley.

[†]Supported in part by the Belgium National Fund for Scientific Research and GOA-grant GOA/2007/04 of Research Fund K.U.Leuven.

[‡]Supported in part by the National Institutes of Health.

[§]Supported in part by the National Institutes of Health and the National Science Foundation.

[1]Department of Mathematics and Center for Statistics, Katholieke Universiteit Leuven, W. de Croylaan 54, B-3001 Leuven (Heverlee), Belgium, e-mail: **irene.gijbels@wis.kuleuven.be**

[2]Department of Biostatistics, CB 7420, University of North Carolina, Chapel Hill, NC 27599-7420, USA, e-mail: **lin@bios.unc.edu**

[3]Department of Statistics, Columbia University, 1255 Amsterdam Avenue, 10th Floor, New York, NY 10027, USA, e-mail: **zying@stat.columbia.edu**

or by solving $\{S(\beta) = 0\}$, where

$$(1.2) \qquad S(\beta) = \sum_{i=1}^{n} \int_0^{\infty} \left\{ Z_i(t) - \frac{\sum_{j=1}^{n} 1_{(X_j \geq t)} e^{\beta' Z_j(t)} Z_j(t)}{\sum_{j=1}^{n} 1_{(X_j \geq t)} e^{\beta' Z_j(t)}} \right\} \delta_i d1_{(X_i \leq t)}.$$

Under suitable regularity conditions, $n^{-1/2} S(\beta_0) \xrightarrow{d} \mathcal{N}(0, V)$ and $n^{1/2}(\hat{\beta}_f - \beta_0) \xrightarrow{d} \mathcal{N}(0, V^{-1})$, where $V = -\lim_{n \to \infty} n^{-1} \partial S(\beta_0)/\partial \beta$ [1]. These asymptotic properties provide the basis for making inference about $\beta_0$. For the one-dimensional (dichotomous) $Z$, the nonparametric test based on $S(0)$ for testing $\beta_0 = 0$ has been better known as the (two-sample) log rank test.

The estimation of the cumulative hazard function $\Lambda(t) = \int_0^t \lambda(s)ds$ and the survival function $\overline{F}(t) = e^{-\Lambda(t)}$ is also of interest. In the one-sample case, where no covariates are modeled, $\Lambda(t)$ is commonly estimated by the Nelson-Aalen estimator

$$(1.3) \qquad \hat{\Lambda}_{NA}(t) = \int_0^t \frac{\sum_{i=1}^{n} \delta_i d1_{(X_i \leq s)}}{\sum_{j=1}^{n} 1_{(X_j \geq s)}},$$

and the corresponding survival function estimator $\hat{\overline{F}}_{NA}(t) = e^{-\hat{\Lambda}_{NA}(t)}$ is asymptotically equivalent to the well-known Kaplan-Meier estimator

$$(1.4) \qquad \hat{\overline{F}}_{KM}(t) = \prod_{X_i \leq t} \left\{ 1 - \frac{\delta_i}{\sum_{j=1}^{n} 1_{(X_j \geq X_i)}} \right\}.$$

Motivated by the Nelson-Aalen estimator, Breslow [2] suggested that the cumulative baseline hazard function $\Lambda_0(t) = \int_0^t \lambda_0(s)ds$ under the Cox model be estimated by

$$(1.5) \qquad \hat{\Lambda}_B(t) = \int_0^t \frac{\sum_{i=1}^{n} \delta_i d1_{(X_i \leq s)}}{\sum_{j=1}^{n} 1_{(X_j \geq s)} e^{\hat{\beta}_f' Z_j(s)}}.$$

Both $n^{1/2}\{\hat{\Lambda}_{NA}(\cdot) - \Lambda(\cdot)\}$ and $n^{1/2}\{\hat{\Lambda}_B(\cdot) - \Lambda_0(\cdot)\}$ converge weakly to zero-mean Gaussian processes [1, 3, 6, 14].

All of the aforementioned procedures assume complete measurements on the failure indicators $\delta_i$ $(i = 1, \ldots, n)$. In many applications, however, the values of $\{\delta_i\}$ are missing for some study subjects. We shall distinguish between two types of missingness. For Type I missingness, $\{\delta_i\}$ are missing completely at random among all subjects. For Type II missingness, $\{\delta_i\}$ take value 0 for some subjects and are missing completely at random among the remaining subjects. By missing completely at random, we mean that the missing mechanism is independent of everything else. The following two examples demonstrate how such missingness arises in practice.

**Example 1.** (Type I missingness). Suppose that a series system has two independent components I and II and let $T$ and $C$ represent times to failure of I and II respectively. The potential observations for a single system consist of $X = T \wedge C$ and $\delta = 1_{(T \leq C)}$. Suppose that a large number of systems are operated until failure. Also suppose that the diagnosis of a system to identify which component failed is so costly that it can only be done for a random sample of the systems under testing. Thus we observe all $\{X_i\}$ and a random subset of $\{\delta_i\}$.

**Example 2.** (Type II missingness). In the medical study, investigators are often interested in the time to death attributable to a particular disease, in which case $\delta_i = 1$ if and only if the $i$th subject died from that disease. Typically, the causes of death are unknown for some deaths because it requires extra efforts (e.g., performing autopsies or obtaining death certificates) to gather such information. Thus the values of $\{\delta_i\}$ may be missing among the deaths. On the other hand, if the $i$th subject has been withdrawn from the study before its termination or is still alive at the end of the study, then $\delta_i$ must be 0. Hence, we have Type II missingness provided that the deaths with known causes are representative of all the subjects who died.

The most commonly adopted strategy for handling missing values is the complete-case analysis, which totally disregards all the subjects with unknown failure statuses. This approach is valid under Type I missingness; however, it can be highly inefficient if there is heavy missingness. For Type II missingness, the complete-case analysis does not even yield consistent estimators.

There have been a few articles on estimating the survival distribution of a homogeneous population in the presence of missing failure indicators. Notably, [5] used the nonparametric maximum likelihood method in conjunction with the EM algorithm to derive an estimator that is analogous to the Kaplan-Meier estimator (1.4). According to [10], however, the maximum likelihood as well as the self-consistent estimators are in general nonunique and inconsistent. Two alternative estimators are proposed in [10] under Type I missingness. As will be discussed in Section 3, these estimators have some undesirable properties. On the more challenging regression problem, there has been little progress. The only solution seems to have been the modified log rank test for Example 1.2 proposed [8]. As admitted by these authors, they made some rather unrealistic assumptions, including the independence between the covariate and the causes of death not under study as well as the proportionality of the hazard rate for the cause of interest and that of the other causes. On the other hand, further developments along the line of efficient estimation can be found in [11, 13, 15]. Furthermore, [17] deals with the additive hazards regression model.

This paper provides a treatment of the Cox regression analysis and the survival function estimation under both types of missingness. In the next section, we introduce a class of estimating functions for $\beta_0$ under Type I missingness which incorporates the partial information from the individuals with unknown $\delta_i$. The consistency and asymptotic normality of the resulting estimators are established. A simple adaptive estimator is constructed which has the smallest variance-covariance matrix among the proposed class of estimators including the complete-case estimator. Simulation studies show that the adaptive estimator is suitable for practical use. Section 3 deals with the survival function estimation under Type I missingness. For the one-sample case, we derive an adaptive estimator which offers considerable improvements over the complete-case and Lo's estimators [10]. Estimation of the cumulative baseline hazard function for the Cox model is also studied. In Section 4, we apply the ideas developed in Sections 2 and 3 to Type II missingness to obtain consistent estimators with similar optimality properties. Note that some of the technical developments there are streamlined and may be traced to a technical report [7]. We conclude this paper with some discussions in Section 5.

## 2. Cox regression under Type I missingness

In this section, we propose estimating functions for the parameter vector $\beta_0$ which utilize the partial information from the subjects with unknown failure indicators.

The asymptotic properties of these functions and the resulting parameter estimators are studied in detail. Throughout the paper, we shall make the following assumption, which is satisfied in virtually all practical situations.

*Boundedness condition.* The covariate processes $Z_i(\cdot) = \{Z_{i1}(\cdot), \ldots, Z_{ip}(\cdot)\}'$ ($i = 1, \ldots, n$) are of bounded variation with a uniform bound, i.e., there exists $K > 0$ such that for all $i$,

$$\sum_{j=1}^{p} \left\{ |Z_{ij}(0)| + \int_0^\infty |dZ_{ij}(t)| \right\} \le K.$$

Let $\xi_i$ indicate, by the value 1 vs. 0, whether $\delta_i$ is known or not. Under Type I missingness, the data consist of i.i.d random vectors $(X_i, \xi_i, \xi_i\delta_i, Z_i')$ ($i = 1, \ldots, n$), where $\xi_i$ is independent of $(X_i, \delta_i, Z_i')$ for every $i$. Write $\rho = P(\xi_1 = 1)$.

Note that the partial likelihood score function (1.2) is the sum over all the observed failure times of the differences between the covariate vectors of the subjects who fail and the weighted averages of the covariate vectors among the subjects under observation. In view of this fact, we introduce the following estimating function:

$$(2.1) \qquad S_1(\beta) = \sum_{i=1}^{n} \int_0^\infty \left\{ Z_i(t) - \bar{Z}(\beta, t) \right\} \xi_i dN_i^u(t),$$

where $\bar{Z}(\beta, t) = \sum_{j=1}^{n} 1_{(X_j \ge t)} e^{\beta' Z_j(t)} Z_j(t) / \sum_{j=1}^{n} 1_{(X_j \ge t)} e^{\beta' Z_j(t)}$ and $N_i^u(t) = \delta_i 1_{(X_i \le t)}$. In the sequel, we shall also use the notation $Y_i(t) = 1_{(X_i \ge t)}$, $N_i(t) = 1_{(X_i \le t)}$ and $N_i^c(t) = (1 - \delta_i) 1_{(X_i \le t)}$. Note that $\{N_i^u, N_i^c\}$ may not be fully observable whereas $\{N_i, \xi_i N_i^u, \xi_i N_i^c\}$ are always observed. Another way of deriving (2.1) is to modify the partial likelihood function (1.1) by omitting the factors for which the $\delta_i$ are missing. Then $S_1(\beta)$ can be obtained by the usual way of differentiating the "log-likelihood function".

**Theorem 2.1.** *Let* $S_1(\beta, t) = \sum_{i=1}^{n} \int_0^t \left\{ Z_i(s) - \bar{Z}(\beta, s) \right\} \xi_i dN_i^u(s)$.

(i) *The process* $n^{-1/2} S_1(\beta_0, \cdot)$ *converges weakly to a zero-mean Gaussian martingale with variance function*

$$(2.2) \qquad V_1(t) = E\left[ \int_0^t \{Z_1(s) - \bar{z}(\beta_0, s)\}^{\otimes 2} \xi_1 dN_1^u(s) \right],$$

*where* $\bar{z}(\beta, t) = E\left\{ Y_1(t) e^{\beta' Z_1(t)} Z_1(t) \right\} / E\left\{ Y_1(t) e^{\beta' Z_1(t)} \right\}$.

(ii) *Define* $\tilde{\beta}$ *as the root of* $\{S_1(\beta) = 0\}$. *If* $V_1 = V_1(\infty)$ *is nonsingular, then* $n^{1/2}(\tilde{\beta} - \beta_0) \xrightarrow{d} \mathcal{N}(0, V_1^{-1})$.

**Remarks.** (1) It is simple to show that $V_1 = \rho V$, where $V$ is the limiting covariance matrix for $\hat{\beta}_f$ defined in Section 1. By the arguments of [1], $V_1(t)$ can be consistently estimated by

$$\hat{V}_1(t) = n^{-1} \sum_{i=1}^{n} \int_0^t \left\{ \frac{\sum_{j=1}^{n} Y_j(s) e^{\tilde{\beta}' Z_j(s)} Z_j^{\otimes 2}(s)}{\sum_{j=1}^{n} Y_j(s) e^{\tilde{\beta}' Z_j(s)}} - \bar{Z}^{\otimes 2}(\tilde{\beta}, s) \right\} \xi_i dN_i^u(s).$$

(2) The nonsingularity of $V_1$ is a very mild assumption and it is true in practically all meaningful situations.

(3) The difference between the process $S_1(\beta, t)$ and the partial likelihood score process under the complete-case analysis

$$S_d(\beta, t) = \sum_{i=1}^n \int_0^t \left\{ Z_i(s) - \bar{Z}_d(\beta, s) \right\} \xi_i dN_i^u(s),$$

where $\bar{Z}_d(\beta, t) = \sum_{j=1}^n \xi_j Y_j(t) e^{\beta' Z_j(t)} Z_j(t) / \sum_{j=1}^n \xi_j Y_j(t) e^{\beta' Z_j(t)}$, is that the subjects with unknown failure indicators are included in the calculation of $\bar{Z}$, but not in that of $\bar{Z}_d$. It is somewhat surprising to note that $S_d(\beta, \cdot)$ and the corresponding estimator $\hat{\beta}_d$ have the same asymptotic distributions as those of $S_1(\beta, \cdot)$ and $\tilde{\beta}$, respectively, even though $\bar{Z}(\beta, t)$ is a more accurate estimator of $\bar{z}(\beta, t)$ than $\bar{Z}_d(\beta, t)$ is. As will be seen in the proof of Theorem 2.1, however, $S_d(\beta, \cdot)$ and $\hat{\beta}_d$ themselves are not asymptotically equivalent to $S_1(\beta, \cdot)$ and $\tilde{\beta}$. Simulation results to be reported later in the section reveal that $\tilde{\beta}$ tends to be slightly more efficient than $\hat{\beta}_d$ for small and moderate-sized samples.

(4) The use of $S_1(\beta)$ may incur substantial loss of information, especially when $\rho$ is small, since the asymptotic distribution of $\tilde{\beta}$ is the same as that of $\hat{\beta}_d$, which only uses data with known failure indicators. Indeed, the purpose of this section is to construct a new estimator that combines $S_1(\beta)$ with an estimating function utilizing the counting processes $N_i(\cdot)$ associated with $\xi_i = 0$. In this connection, the estimating function $S_1$ plays only a transitional role.

*Proof of Theorem 2.1.* For notational simplicity, assume $p = 1$. Let $M_i(t) = N_i^u(t) - \int_0^t Y_i(s) e^{\beta_0 Z_i(s)} \lambda_0(s) ds$, which are martingale processes with respect to an appropriate $\sigma$-filtration [1]. Decompose $S_1(\beta_0, t)$ into two parts

$$\begin{aligned} S_1(\beta_0, t) &= \sum_{i=1}^n \int_0^t \left\{ Z_i(s) - \bar{Z}(\beta_0, s) \right\} \xi_i dM_i(s) \\ &\quad + \sum_{i=1}^n (\xi_i - \rho) \int_0^t \left\{ Z_i(s) - \bar{Z}(\beta_0, s) \right\} e^{\beta_0 Z_i(s)} Y_i(s) \lambda_0(s) ds \\ &= S_{11}(t) + S_{12}(t), \qquad \text{say.} \end{aligned}$$

Now $n^{-1/2} S_{11}(\cdot)$ is a martingale. By the arguments of [1], $n^{-1/2} S_{11}(t)$ is asymptotically equivalent to $n^{-1/2} \tilde{S}_{11}(t) = n^{-1/2} \sum_{i=1}^n \int_0^t W_i(s) \xi_i dM_i(s)$, where $W_i(s) = Z_i(s) - \bar{z}(\beta_0, s)$, and converges weakly in $\mathcal{D}[0, \infty)$ to a Gaussian martingale with variance function $V_1(t)$. Note that the tightness of $n^{-1/2} S_{11}(\cdot)$ at $\infty$ can be easily handled along the lines of [6]. From Lemma 1(i) given at the end of the section, $n^{-1/2} S_{12}(t)$ is also tight and is asymptotically equivalent to

$$n^{-1/2} \tilde{S}_{12}(t) = n^{-1/2} \sum_{i=1}^n (\xi_i - \rho) \int_0^t W_i(s) e^{\beta_0 Z_i(s)} Y_i(s) \lambda_0(s) ds.$$

Hence, $n^{-1/2} S_1(\beta_0, \cdot)$ is asymptotically equivalent to $n^{-1/2} \left\{ \tilde{S}_{11}(\cdot) + \tilde{S}_{12}(\cdot) \right\}$, which converges weakly to a zero-mean Gaussian process with covariance function at $(t, t')$ that can be shown to be equal to

$$n^{-1} E \left[ \left\{ \tilde{S}_{11}(t) + \tilde{S}_{12}(t) \right\} \left\{ \tilde{S}_{11}(t') + \tilde{S}_{12}(t') \right\} \right] = V_1(t \wedge t').$$

To prove part (ii) of the theorem, note that $-n^{-1} \partial S_1(\beta) / \partial \beta$ is positive (positive definite for $p > 1$) and converges to $E \left[ \int_0^\infty \left\{ Z_1(t) - \bar{z}(\beta, t) \right\}^2 \xi_1 dN_1^u(t) \right]$. Thus, $\tilde{\beta}$ is uniquely defined and the arguments of [1] entail the convergence of $n^{1/2}(\tilde{\beta} - \beta_0)$ to $\mathcal{N}(0, V_1^{-1})$.                                                                                  □

To incorporate the partial survival information from those subjects with missing $\delta_i$, it is natural to consider the counting processes $(1 - \xi_i)N_i(\cdot)$ and to subtract off the jumps due to censoring. In this connection, we introduce

(2.3)
$$S_2(\beta, t) = \sum_{i=1}^{n} \int_0^t \left\{ Z_i(s) - \bar{Z}(\beta, s) \right\} \left\{ (1 - \xi_i) dN_i(s) - \hat{\rho}^{-1}(1 - \hat{\rho})\xi_i dN_i^c(s) \right\},$$

where $\hat{\rho} = n^{-1}\sum_{i=1}^{n}\xi_i$, noting that $E\left\{(1 - \xi_1)N_1(t) - \rho^{-1}(1 - \rho)\xi_1 N_1^c(t)\right\} = E\{N_1^u(t)\}$.

**Theorem 2.2.** *The process $n^{-1/2}S_2(\beta_0, \cdot)$ is asymptotically independent of $n^{-1/2} \times S_1(\beta_0, \cdot)$ and converges weakly to a zero-mean Gaussian process with covariance function*

$$V_2(t, t') = E\left\{ \int_0^{t \wedge t'} W_1^{\otimes 2}(s)(1 - \xi_1)dN_1^u(s) \right\}$$
$$+ \rho^{-1}(1 - \rho)E\left[ \left\{ N_1^{CZ}(t) - EN_1^{CZ}(t) \right\} \left\{ N_1^{CZ}(t') - EN_1^{CZ}(t') \right\}' \right],$$

*where $N_i^{CZ}(t) = \int_0^t \left\{ Z_i(s) - \bar{z}(\beta_0, s) \right\} dN_i^c(s)$ $(i = 1, \ldots, n)$.*

*Proof.* Again assume $p = 1$. Since $\hat{\rho} - \rho = O_p(n^{-1/2})$, by the usual delta method,

$$S_2(\beta_0, t) = \left[ \sum_{i=1}^{n} \int_0^t \left\{ Z_i(s) - \bar{Z}(\beta_0, s) \right\} (1 - \xi_i)dM_i(s) \right.$$
$$\left. + \sum_{i=1}^{n} \int_0^t \left\{ Z_i(s) - \bar{Z}(\beta_0, s) \right\} \left\{ 1 - \xi_i - (1 - \rho) \right\} e^{\beta_0 Z_i(s)} Y_i(s)\lambda_0(s)ds \right]$$
$$+ \left[ \sum_{i=1}^{n} \int_0^t \left\{ Z_i(s) - \bar{Z}(\beta_0, s) \right\} dN_i^c(s)\rho^{-1}(\rho - \xi_i) \right.$$
$$\left. + \sum_{i=1}^{n} \int_0^t \left\{ Z_i(s) - \bar{Z}(\beta_0, s) \right\} dN_i^c(s)\rho^{-1}(\hat{\rho} - \rho) \right]$$
$$+ r_n(t)$$
$$= S_{21}(t) + S_{22}(t) + r_n(t), \qquad \text{say.}$$

Here the remainder term $r_n$ is uniformly negligible in the sense that $\sup_t |r_n(t)| = o_p(n^{1/2})$. Note that $S_{21}(t)$ is the same as $S_1(t)$ except that $\{\xi_i\}$ there are replaced by $\{1 - \xi_i\}$. Thus $S_{21}(t)$ is tight in $\mathcal{D}[0, \infty)$ and is asymptotically equivalent to

$$\tilde{S}_{21}(t) = \sum_{i=1}^{n} \int_0^t W_i(s)(1 - \xi_i)dM_i(s) + \sum_{i=1}^{n} \int_0^t W_i(s)(\rho - \xi_i)e^{\beta_0 Z_i(s)}Y_i(s)\lambda_0(s)ds.$$

By Lemma 1(ii), $S_{22}(t)$ is asymptotically equivalent to

$$\tilde{S}_{22}(t) = -\sum_{i=1}^{n} \rho^{-1}(\xi_i - \rho) \left\{ N_i^{CZ}(t) - EN_i^{CZ}(t) \right\}.$$

By writing

$$\tilde{S}_{21}(t) = \sum_{i=1}^{n} \int_0^t W_i(s)(1-\rho)dM_i(s) + \sum_{i=1}^{n} \int_0^t W_i(s)dN_i^u(s)(\rho - \xi_i),$$

we can show that, for any $t$ and $t'$,

$$(2.4) \qquad\qquad E\left\{ \tilde{S}_{21}(t)\tilde{S}_{22}(t') \right\} = 0.$$

Thus $n^{-1/2}\left\{ \tilde{S}_{21}(\cdot) + \tilde{S}_{22}(\cdot) \right\}$ converges weakly to a zero-mean Gaussian process with $V_2$ as its covariance function.

Similar to (2.4), $E\left[ \{\tilde{S}_{11}(t) + \tilde{S}_{12}(t)\}\tilde{S}_{22}(t') \right] = 0$ for any $t$ and $t'$. Thus to prove the asymptotic independence between $S_1$ and $S_2$, it suffices to show

$$E\left[ \{\tilde{S}_{11}(t) + \tilde{S}_{12}(t)\}\tilde{S}_{21}(t') \right] = 0.$$

To this end, we can apply the same covariance calculation as employed in the proof of Theorem 2.1 to show that

$$\begin{aligned}
&E\left[ \{\tilde{S}_{11}(t) + \tilde{S}_{12}(t)\}\tilde{S}_{21}(t') \right]\\
&\quad = nE\left\{ \int_0^t W_1(s)\xi_1 dM_1(s) \int_0^{t'} W_1(s')(1-\xi_1)dM_1(s') \right\} = 0.
\end{aligned}$$

$\square$

By combining $S_1$ and $S_2$, more efficient estimators of $\beta_0$ may be obtained. Specifically, given a $p \times p$ matrix $D$, we can define $\hat{\beta}$ as a solution to

$$(2.5) \qquad\qquad S_1(\beta) + DS_2(\beta) = 0.$$

**Theorem 2.3.** *Suppose that $\{\rho V + (1-\rho)DV\}$ is nonsingular. Let $V_2 = V_2(\infty, \infty)$. Then $n^{1/2}(\hat{\beta} - \beta_0) \xrightarrow{d} \mathcal{N}(0, \Sigma(D))$, where*

$$(2.6) \qquad \Sigma(D) = \{\rho V + (1-\rho)DV\}^{-1} (\rho V + DV_2 D') \{\rho V + (1-\rho)VD'\}^{-1}.$$

*In particular, $D^* = (1-\rho)VV_2^{-1}$ yields*

$$\Sigma(D^*) = \left\{ \rho V + (1-\rho)^2 VV_2^{-1}V \right\}^{-1}$$

*and is optimal in the sense that $\Sigma(D) - \Sigma(D^*)$ is nonnegative definite for any $D$.*

**Remarks.** (1) Let $V_{CZ} = E\left\{ N_1^{CZ}(\infty) - EN_1^{CZ}(\infty) \right\}^{\otimes 2}$. Then $V_2 = (1-\rho)V + \rho^{-1}(1-\rho)V_{CZ}$. For $p = 1$, $D^* = V/(V + \rho^{-1}V_{CZ})$ and $\Sigma(D^*) = (V + \rho^{-1}V_{CZ})/\{V(V + V_{CZ})\}$. This variance will be close to the ideal $V^{-1}$ if either $\rho$ is close to 1 (light missingness) or $V_{CZ}$ is close to zero (light censorship).

(2) A consistent estimator for $\Sigma(D)$ may be obtained by replacing $\rho$, $V$ and $V_2$ in (2.6) by $\hat{\rho}$, $\hat{V}(\hat{\beta})$ and $\hat{V}_2(\hat{\beta})$, where

$$\hat{V}(\beta) = \frac{1}{\sum_{i=1}^{n} \xi_i} \sum_{i=1}^{n} \int_0^\infty \left\{ \frac{\sum_{j=1}^{n} Y_j(t)e^{\beta' Z_j(t)}Z_j^{\otimes 2}(t)}{\sum_{j=1}^{n} Y_j(t)e^{\beta' Z_j(t)}} - \bar{Z}^{\otimes 2}(\beta, t) \right\} \xi_i dN_i^u(t),$$

*Monte Carlo estimates for the sampling means and variances of four estimators of $\beta_0$ and for the sizes of the corresponding 0.05-level Wald tests for testing $H_0 : \beta_0 = 0$ under the Model $\lambda(t|Z) = 1$*

| $\rho$ | Estimator | 20% Censoring | | | 50% Censoring | | | 70% Censoring | | |
|---|---|---|---|---|---|---|---|---|---|---|
| | | Mean | Var. | Size | Mean | Var. | Size | Mean | Var. | Size |
| 0.8 | $\hat{\beta}_f$ | –0.001 | 0.015 | 0.056 | –0.001 | 0.023 | 0.054 | 0.002 | 0.038 | 0.052 |
| | $\hat{\beta}_d$ | –0.001 | 0.018 | 0.053 | –0.002 | 0.030 | 0.052 | 0.002 | 0.049 | 0.049 |
| | $\tilde{\beta}$ | –0.001 | 0.018 | 0.054 | –0.001 | 0.029 | 0.053 | 0.002 | 0.049 | 0.050 |
| | $\hat{\beta}^*$ | –0.001 | 0.015 | 0.056 | –0.002 | 0.027 | 0.055 | 0.002 | 0.046 | 0.051 |
| 0.5 | $\hat{\beta}_f$ | –0.001 | 0.015 | 0.056 | –0.001 | 0.023 | 0.054 | 0.002 | 0.038 | 0.052 |
| | $\hat{\beta}_d$ | –0.002 | 0.032 | 0.057 | –0.002 | 0.054 | 0.053 | 0.001 | 0.092 | 0.049 |
| | $\tilde{\beta}$ | 0.001 | 0.029 | 0.050 | –0.001 | 0.048 | 0.052 | 0.003 | 0.082 | 0.050 |
| | $\hat{\beta}^*$ | –0.001 | 0.017 | 0.056 | –0.002 | 0.037 | 0.052 | 0.002 | 0.071 | 0.048 |

NOTE: $Z$ is standard normal. The censoring time is exponentially distributed with hazard rate $\lambda^c$, where $\lambda^c$ is chosen to achieve the desired censoring percentage. The sample size is 100. Each block is based on 10,000 replications. The random number generator of [16] is used.

$$\hat{V}_2(\beta) = (1 - \hat{\rho})\hat{V}(\beta) + \hat{\rho}^{-1}(1 - \hat{\rho})\hat{V}_{CZ}(\beta),$$

$$\hat{V}_{CZ}(\beta) = \frac{1}{\sum_{i=1}^{n} \xi_i} \sum_{i=1}^{n} \int_0^\infty \left\{ Z_i(t) - \bar{Z}(\beta, t) \right\}^{\otimes 2} \xi_i dN_i^c(t)$$

$$- \left[ \frac{1}{\sum_{i=1}^{n} \xi_i} \sum_{i=1}^{n} \int_0^\infty \left\{ Z_i(t) - \bar{Z}(\beta, t) \right\} \xi_i dN_i^c(t) \right]^{\otimes 2}.$$

Since we can estimate the optimal weight $D^*$ consistently by $\hat{D}^* = (1 - \hat{\rho})\hat{V}(\tilde{\beta}) \times \hat{V}_2^{-1}(\tilde{\beta})$, an "adaptive" estimator of $\beta_0$ that achieves the lower variance-covariance bound $\Sigma(D^*)$ may be constructed. Specifically, we can first use $\tilde{\beta}$ from $\{S_1(\beta) = 0\}$ to compute $\hat{D}^*$ and then obtain the adaptive estimator by solving

$$(2.7) \qquad\qquad S_1(\beta) + \hat{D}^* S_2(\beta) = 0.$$

**Corollary 1.** *Let $\hat{\beta}^*$ be the estimator given by (2.7). Then under the same assumptions as Theorem 2.3, $n^{1/2}(\hat{\beta}^* - \beta_0) \xrightarrow{d} \mathcal{N}(0, \Sigma(D^*))$. In addition, $\Sigma(D^*)$ can be consistently estimated by $\left\{ \hat{\rho}\hat{V}(\hat{\beta}^*) + (1 - \hat{\rho})^2 \hat{V}(\hat{\beta}^*) \hat{V}_2^{-1}(\hat{\beta}^*) \hat{V}(\hat{\beta}^*) \right\}^{-1}$.*

We have carried out extensive Monte Carlo experiments to investigate the finite-sample behaviour of the proposed adaptive estimator $\hat{\beta}^*$ and to compare it with the full-data estimator $\hat{\beta}_f$, the complete-case estimator $\hat{\beta}_d$ and the $S_1(\beta)$ estimator $\tilde{\beta}$. The key results are summarized in Tables 1 and 2. The biases of all four estimators and of their variance estimators (the latter not shown here) are negligible, and the associated Wald tests have proper sizes. The adaptive estimator is always more efficient than $\hat{\beta}_d$ and $\tilde{\beta}$, as is reflected in the sampling variances of the estimators as well as in the powers of the Wald tests. The gains in the relative efficiencies increase as the missing probability increases and decrease as the censoring probability increases. The efficiency of $\hat{\beta}^*$ relative to $\hat{\beta}_f$ is close to 1 when censoring is light. The estimator $\tilde{\beta}$ seems to have slightly better small-sample efficiency than $\hat{\beta}_d$.

*Proof of Theorem 2.3.* From its definition, $-\partial S_1(\beta)/\partial \beta$ is, with probability 1, positive definite. Thus, following [1], we can show that $(\rho n)^{-1} S_1(\beta)$ converges uniformly

TABLE 2

*Monte Carlo estimates for the sampling means and variances of four estimators of $\beta_0$ and for the powers of the corresponding 0.05-level wald tests for testing $h_0 : \beta_0 = 0$ under the Model $\lambda(t|Z) = e^{0.5Z}$*

| $\rho$ | Estimator | 20% Censoring | | | 50% Censoring | | | 70% Censoring | | |
|---|---|---|---|---|---|---|---|---|---|---|
| | | Mean | Var. | Power | Mean | Var. | Power | Mean | Var. | Power |
| 0.8 | $\hat{\beta}_f$ | 0.509 | 0.017 | 0.984 | 0.511 | 0.026 | 0.912 | 0.514 | 0.041 | 0.755 |
| | $\hat{\beta}_d$ | 0.511 | 0.021 | 0.956 | 0.514 | 0.033 | 0.844 | 0.518 | 0.053 | 0.655 |
| | $\tilde{\beta}$ | 0.510 | 0.021 | 0.955 | 0.513 | 0.033 | 0.844 | 0.516 | 0.052 | 0.653 |
| | $\hat{\beta}^*$ | 0.509 | 0.018 | 0.980 | 0.512 | 0.030 | 0.876 | 0.516 | 0.049 | 0.684 |
| 0.5 | $\hat{\beta}_f$ | 0.509 | 0.017 | 0.984 | 0.511 | 0.026 | 0.912 | 0.514 | 0.041 | 0.755 |
| | $\hat{\beta}_d$ | 0.516 | 0.038 | 0.813 | 0.522 | 0.061 | 0.624 | 0.530 | 0.102 | 0.435 |
| | $\tilde{\beta}$ | 0.511 | 0.035 | 0.821 | 0.514 | 0.054 | 0.642 | 0.520 | 0.088 | 0.450 |
| | $\hat{\beta}^*$ | 0.509 | 0.021 | 0.960 | 0.512 | 0.042 | 0.757 | 0.518 | 0.077 | 0.514 |

NOTE: See NOTE of Table 1.

in any compact set to the nonrandom function

$$m(\beta) = E\left[\int_0^\infty \left\{Z_1(t) - \bar{z}(\beta, t)\right\} dN_1^u(t)\right].$$

For $S_2(\beta)$, we have, by the law of large numbers,

$$n^{-1} \sum_{i=1}^n \int_0^\infty \left\{Z_i(t) - \bar{Z}(\beta, t)\right\} (1 - \xi_i) dN_i(t)$$

$$- n^{-1} \sum_{i=1}^n \int_0^\infty \left\{Z_i(t) - \bar{Z}(\beta, t)\right\} (1 - \rho) dN_i(t)$$

$$= o_p(1) + \int_0^\infty \bar{Z}(\beta, t) d\left\{n^{-1} \sum_{i=1}^n N_i(t)(\xi_i - \rho)\right\}$$

$$= o_p(1),$$

where the last equality follows from the facts that $\sup_t \left|n^{-1} \sum_{i=1}^n N_i(t)(\xi_i - \rho)\right| = o_p(n^{-1/4})$ and that the total variation of $\bar{Z}(\beta, \cdot)$ is at most $O(\log n)$ uniformly for $\beta$ in any compact region. Thus the order $o_p(1)$ is also uniform. Continuing this line of arguments, we get

$$n^{-1} S_2(\beta) = n^{-1} \sum_{i=1}^n \int_0^\infty \left\{Z_i - \bar{Z}(\beta, t)\right\} (1 - \rho) dN_i^u(t) + o_p(1)$$

$$= (1 - \rho) m(\beta) + o_p(1)$$

with the same uniformity. Thus $n^{-1} \left\{S_1(\beta) + DS_2(\beta)\right\}$ is uniformly approximated by $\{\rho I + (1 - \rho)D\} m(\beta)$, which has a unique root $\beta_0$. Hence, $\hat{\beta} \xrightarrow{p} \beta_0$.

The asymptotic normality is easier to show now. Taking the Taylor series expansion of $\left\{S_1(\hat{\beta}) + DS_2(\hat{\beta})\right\}$ at $\beta_0$, we get

$$n^{1/2}(\hat{\beta} - \beta_0) = \left[\left\{\rho I + (1 - \rho)D\right\} V\right]^{-1} n^{-1/2} \left\{S_1(\beta_0) + DS_2(\beta_0)\right\} + o_p(1),$$

which, by Theorems 2.1 and 2.2 and a straightforward matrix manipulation, converges to the desired normal distribution.

To verify the optimality of $D^*$, we note that the estimating function can be linearized around $\beta_0$ and the limiting normal random vectors may be used in place of $n^{-1/2}S_k(\beta_0)$ $(k = 1, 2)$. Specifically, we can consider the following "limiting" linear model

$$\begin{cases} S_1^* = \rho V b + S_{01}^*, \\[2mm] S_2^* = (1 - \rho)V b + S_{02}^*, \end{cases}$$

where $S_{0k}^*$ $(k = 1, 2)$ are independent $\mathcal{N}(0, V_k)$ $(k = 1, 2)$. Recall that $V_1 = \rho V$ and $V_2 = (1 - \rho)(V + \rho^{-1}V_{CZ})$. By the Gauss-Markov theorem, the best linear estimator is

$$\hat{b}^* = \left\{ \rho V + (1 - \rho)^2 V V_2^{-1} V \right\}^{-1} \left\{ S_1^* + (1 - \rho) V V_2^{-1} S_2^* \right\}$$

with variance-covariance matrix $\left\{ \rho V + (1 - \rho)^2 V V_2^{-1} V \right\}^{-1}$, which is exactly $\Sigma(D^*)$.  $\square$

**Lemma 1.** (i) *The process*

$$n^{-1/2}S_{12}(t) = n^{-1/2} \sum_{i=1}^n (\xi_i - \rho) \int_0^t \left\{ Z_i(s) - \bar{Z}(\beta_0, s) \right\} e^{\beta_0' Z_i(s)} \lambda_0(s) ds$$

*is tight in* $\mathcal{D}[0, \infty)$ *and is asymptotically equivalent to*

$$n^{-1/2}\tilde{S}_{12}(t) = n^{-1/2} \sum_{i=1}^n (\xi_i - \rho) \int_0^t \left\{ Z_i(s) - \bar{z}(\beta_0, s) \right\} e^{\beta_0' Z_i(s)} \lambda_0(s) ds$$

*in the sense that* $\sup_t n^{-1/2} \| \tilde{S}_{12}(t) - S_{12}(t) \| = o_p(1)$.

(ii) *The process* $n^{-1/2} \sum_{i=1}^n (\xi_i - \rho) \int_0^t \left\{ Z_i(s) - \bar{Z}(\beta_0, s) \right\} dN_i^c(s)$ *is tight and asymptotically equivalent to*

$$n^{-1/2} \sum_{i=1}^n (\xi_i - \rho) \int_0^t \left\{ Z_i(s) - \bar{z}(\beta_0, s) \right\} dN_i^c(s).$$

*Proof.* Without loss of generality, assume $p = 1$. For any $t_1 < t_2$, we have

$$\begin{aligned}
&\varlimsup_{n \to \infty} E \left\{ n^{-1/2}S_{12}(t_2) - n^{-1/2}S_{12}(t_1) \right\}^4 \\
&= \varlimsup_{n \to \infty} n^{-2} \sum_{i \neq j} E \left[ (\xi_i - \rho) \int_{t_1}^{t_2} \left\{ Z_i(s) - \bar{Z}(\beta_0, s) \right\} e^{\beta_0 Z_i(s)} Y_i(s) \lambda_0(s) ds \right. \\
&\qquad \left. \times (\xi_j - \rho) \int_{t_1}^{t_2} \left\{ Z_j(s) - \bar{Z}(\beta_0, s) \right\} e^{\beta_0 Z_j(s)} Y_j(s) \lambda_0(s) ds \right]^2 \\
&\leq (2K)^4 \varlimsup_{n \to \infty} n^{-2} \sum_{i \neq j} E \left\{ (\xi_i - \rho) \int_{t_1}^{t_2} e^{\beta_0 Z_i(s)} Y_i(s) \lambda_0(s) ds \right. \\
&\qquad \left. \times (\xi_j - \rho) \int_{t_1}^{t_2} e^{\beta_0 Z_j(s)} Y_j(s) \lambda_0(s) ds \right\}^2 \\
&= (2K)^4 \left\{ \rho(1 - \rho) \right\}^2 \left[ E \left\{ \int_{t_1}^{t_2} Y_1(s) e^{\beta_0 Z_1(s)} \lambda_0(s) ds \right\}^2 \right]^2.
\end{aligned}$$

Since $\mu[t_1, t_2] = E\left\{\int_{t_1}^{t_2} Y_1(s)e^{\beta_0 Z_1(s)}\lambda_0(s)ds\right\}^2$ is a finite measure on $[0, \infty)$, the moment criterion ([12], page 52, formula (30)) implies the tightness of $n^{-1/2}S_{12}$. Likewise, $n^{-1/2}\tilde{S}_{12}$ is also tight. Furthermore, let $t_0 = \inf\{t : EY_1(t) = 0\}$. Then it is easy to see that $\sup_{s \leq t}|\bar{Z}(\beta_0, s) - \bar{z}(\beta_0, s)| \xrightarrow{p} 0$ for any $t < t_0$. Thus the equivalence of $n^{-1/2}S_{12}$ and $n^{-1/2}\tilde{S}_{12}$ follows from the tightness just proved. Hence (i) holds.

The proof of (ii) is very much the same as that of (i). Because of possible discontinuity of $EN_1^c(t)$ in $t$, another moment condition ([12], page 51, formula (25)) should be used. Note that the tightness continues to hold even if the measure $\mu$ there is discontinuous. □

## 3. Cumulative hazard function estimation under Type I missingness

In this section, we first deal with the problem of nonparametric estimation of the cumulative hazard function for a homogeneous population under Type I missingness. We shall discuss the estimators proposed in [10] and give our own solutions. We then apply the ideas to the estimation of the cumulative baseline hazard function for the Cox model. In both cases, asymptotic distributions of the relevant estimators are derived.

In the one-sample case, the observations consist of i.i.d. random vectors $(X_i, \xi_i, \xi_i\delta_i)$ $(i = 1, \ldots, n)$, where $X_i = T_i \wedge C_i$, $\delta_i = 1_{(T_i \leq C_i)}$ and $\xi_i$ is the missing indicator independent of $(X_i, \delta_i)$. Assume that $T_i$ is independent of $C_i$ and that $T_i$ has a continuous distribution function. Let $F(t) = P(T_1 \leq t)$, $\Lambda(t) = \int_0^t dF(s)/\{1 - F(s)\}$, $G(t) = P(C_1 \leq t)$, $\Lambda_G(t) = \int_0^t dG(s)/\{1 - G(s-)\}$, $H(t) = \{1 - F(t)\}\{1 - G(t-)\}$, $A(t) = \int_0^t d\Lambda(s)/H(s)$ and $A_G(t) = \int_0^t dG(s)/\{(1 - G(s-))H(s)\}$. The notation for $Y_i, N_i, N_i^u, N_i^c, \rho, \hat{\rho}$, etc. introduced in Section 2 will also be used.

Under the setup described above, [10] shows that the nonparametric maximum likelihood method typically does not yield a consistent estimator for $F$, indicating that this is far more complicated than the complete-data situation. Two alternative estimators, $\hat{F}_1$ and $\hat{F}_B$, are also proposed there. It can be shown, by expanding $\log(1 - \hat{F}_A)$, that $\hat{F}_A$ is not a consistent estimator; in particular, Theorem 3 of [10] is not valid. In our notation, the second estimator is given by

$$(3.1) \qquad \hat{F}_B(t) = 1 - \prod_{X_i \leq t}\left\{1 - \frac{1}{\sum_{j=1}^n Y_j(X_i)}\right\}^{\xi_i\delta_i/\hat{\rho}}.$$

Motivated by (3.1), we modify (1.3) to obtain the following estimator for $\Lambda(t)$:

$$(3.2) \qquad \hat{\Lambda}_1(t) = \int_0^t \frac{\sum_{i=1}^n \xi_i dN_i^u(s)}{\hat{\rho}\sum_{j=1}^n Y_j(s)}.$$

By the exponentiation formula of Doleans-Dade ([1], p. 897), the corresponding estimator for $F(t)$ is

$$(3.3) \qquad \hat{F}_1(t) = 1 - \prod_{X_i \leq t}\left\{1 - \frac{\xi_i\delta_i}{\hat{\rho}\sum_{j=1}^n Y_j(X_i)}\right\}.$$

It is easily seen that $\hat{F}_1$ and $\hat{F}_B$ are asymptotically equivalent; however, the cumulative hazard function approach is more convenient for our later developments.

Expression (3.2) also reveals that $\hat{\Lambda}_1$ (and hence $\hat{F}_B$ and $\hat{F}_1$) does not utilize the counting process information from the subjects with $\xi_i = 0$. To recover this information, we introduce

$$(3.4) \qquad \hat{\Lambda}_2(t) = \int_0^t \frac{\sum_{i=1}^n (1-\xi_i) dN_i(s) - \hat{\rho}^{-1}(1-\hat{\rho}) \sum_{i=1}^n \xi_i dN_i^c(s)}{(1-\hat{\rho}) \sum_{i=1}^n Y_i(s)},$$

which shares the same spirit as estimating function $S_2(\beta)$ given in (2.3). Thus, $\Lambda(t)$ can be estimated by

$$(3.5) \qquad \hat{\Lambda}(\alpha, t) = \alpha \hat{\Lambda}_1(t) + (1-\alpha) \hat{\Lambda}_2(t),$$

where $\alpha \in [0, 1]$.

**Theorem 3.1.** *Let $t_0 < H^{-1}(0)$. Then $n^{1/2} \left\{ \hat{\Lambda}(\alpha, \cdot) - \Lambda(\cdot) \right\}$ converges weakly in $\mathcal{D}[0, t_0]$ to a zero-mean Gaussian process with covariance function*

$$(3.6) \qquad \Gamma_\alpha(t, t') = \frac{\alpha^2}{\rho} \left\{ A(t \wedge t') - (1-\rho) \Lambda(t) \Lambda(t') \right\}$$

$$+ \alpha(1-\alpha) \left\{ 2\Lambda(t)\Lambda(t') + \rho^{-1} \Lambda(t) \Lambda_G(t') + \rho^{-1} \Lambda(t') \Lambda_G(t) \right\}$$

$$+ \frac{(1-\alpha)^2}{1-\rho} \left[ A(t \wedge t') + \rho^{-1} A_G(t \wedge t') \right.$$

$$\left. - \rho \left\{ \Lambda(t) + \rho^{-1} \Lambda_G(t) \right\} \left\{ \Lambda(t') + \rho^{-1} \Lambda_G(t') \right\} \right].$$

*For fixed $t$, $n^{1/2} \left\{ \hat{\Lambda}(\alpha, t) - \Lambda(t) \right\} \xrightarrow{d} \mathcal{N}(0, \Gamma_\alpha(t))$, where*

$$(3.7) \quad \Gamma_\alpha(t) = \frac{\alpha^2}{\rho} \left\{ A(t) - (1-\rho) \Lambda^2(t) \right\} + 2\alpha(1-\alpha) \left\{ \Lambda^2(t) + \rho^{-1} \Lambda(t) \Lambda_G(t) \right\}$$

$$+ \frac{(1-\alpha)^2}{1-\rho} \left[ A(t) + \rho^{-1} A_G(t) - \rho \left\{ \Lambda(t) + \rho^{-1} \Lambda_G(t) \right\}^2 \right],$$

*which reaches its minimum when $\alpha$ equals*

$$\alpha^* = \frac{\rho \left\{ A(t) - \Lambda^2(t) \right\} + A_G(t) - \Lambda_G^2(t) - (1+\rho) \Lambda(t) \Lambda_G(t)}{A(t) - \Lambda^2(t) + A_G(t) - \Lambda_G^2(t) - 2\Lambda(t) \Lambda_G(t)}.$$

**Remarks.** (1) If we choose $\alpha_n \xrightarrow{p} \alpha$, then $\hat{\Lambda}(\alpha_n, \cdot)$ has the same asymptotic distribution as $\hat{\Lambda}(\alpha, \cdot)$. Since $\alpha^*$ can be estimated consistently, the "optimal" estimator of $\Lambda$ can be constructed adaptively. To be specific, $\alpha^*$ may be estimated by

$$\hat{\alpha}^* = \frac{\hat{\rho} \left\{ \hat{A}_1(t) - \hat{\Lambda}_1^2(t) \right\} + \hat{A}_G(t) - \hat{\Lambda}_G^2(t) - (1+\hat{\rho}) \hat{\Lambda}_1(t) \hat{\Lambda}_G(t)}{\hat{A}_1(t) - \hat{\Lambda}_1^2(t) + \hat{A}_G(t) - \hat{\Lambda}_G^2(t) - 2\hat{\Lambda}_1(t) \hat{\Lambda}_G(t)},$$

where $\hat{A}_1(t) = n \int_0^t d\hat{\Lambda}_1(s) / \sum_{j=1}^n Y_j(s)$, and $\hat{\Lambda}_G(t)$ and $\hat{A}_G(t)$ are the obvious analogs of $\hat{\Lambda}_1(t)$ and $\hat{A}_1(t)$.

(2) A consistent estimator for $\Gamma_\alpha(t, t')$ may be obtained by replacing $\rho$, $A$, $\Lambda$, $A_G$ and $\Lambda_G$ in (3.6) by $\hat{\rho}$, $\hat{A}_1$, $\hat{\Lambda}_1$, $\hat{A}_G$ and $\hat{\Lambda}_G$.

(3) Two special cases deserve extra attention. If $\alpha = 1$, then $\hat{\Lambda}(\alpha, t)$ reduces to $\hat{\Lambda}_1(t)$. In that case, the asymptotic variance $\Gamma_1(t) = \rho^{-1} \left\{ A(t) - (1-\rho) \Lambda^2(t) \right\}$, which agrees with Lo's result when the exponentiation is taken into account, and *Simulation summary statistics for the adaptive estimator $\hat{F}(\hat{\alpha}^*, t)$ at $t = F^{-1}(0.5)$ under the exponential model $F(t) = 1 - e^{-t}$*

| | 20% Censoring | | 50% Censoring | | 70% Censoring | |
|---|---|---|---|---|---|---|
| | $\rho = 0.8$ | $\rho = 0.5$ | $\rho = 0.8$ | $\rho = 0.5$ | $\rho = 0.8$ | $\rho = 0.5$ |
| Mean of $\hat{\alpha}^*$ | 0.84 | 0.60 | 0.90 | 0.75 | 0.94 | 0.85 |
| Mean of $\hat{F}(\hat{\alpha}^*, t)$ | 0.497 | 0.497 | 0.496 | 0.495 | 0.493 | 0.490 |
| Var of $n^{1/2}\hat{F}(\hat{\alpha}^*, t)$ | 0.284 | 0.325 | 0.419 | 0.566 | 0.786 | 1.161 |
| Mean of $\hat{V}_{\hat{F}}(\hat{\alpha}^*, t)$ | 0.283 | 0.323 | 0.413 | 0.547 | 0.747 | 1.039 |
| Var of $\hat{F}_d(t)$ / Var of $\hat{F}(\hat{\alpha}^*, t)$ | 1.21 | 1.72 | 1.14 | 1.36 | 1.12 | NA |
| Var of $\hat{F}_B(t)$ / Var of $\hat{F}(\hat{\alpha}^*, t)$ | 1.10 | 1.33 | 1.06 | 1.13 | 1.06 | 1.08 |

NOTE: The censoring time is exponential. The sample size $n = 100$. Each block is based on 10,000 replications. $\hat{V}_{\hat{F}}(\hat{\alpha}^*, t)$ is the variance estimator for $n^{1/2}\hat{F}(\hat{\alpha}^*, t)$, which is $\hat{F}^2(\hat{\alpha}^*, t)$ multiplied by the estimator for $\Gamma_{\alpha^*}(t, t)$ mentioned in Remark (2) of Theorem 3.1. $\hat{F}_d(t)$ is the estimator based on complete cases only and $\hat{F}_B(t)$ is Lo's second estimator. "Mean" and "Var" refer to the sampling mean and variance. NA indicates that the result for the complete-case estimator is not obtainable.

which is less than $\rho^{-1}A(t)$, the variance of the complete-case estimator. On the other hand, if we let $\alpha = \hat{\rho}$, then

$$\hat{\Lambda}(\alpha, t) = \int_0^t \frac{\sum_{i=1}^n \left\{ \xi_i dN_i^u(s) + (1 - \xi_i)dN_i(s) - \hat{\rho}^{-1}(1 - \hat{\rho})\xi_i dN_i^c(s) \right\}}{\sum_{i=1}^n Y_i(s)}$$

with asymptotic variance $\Gamma_\rho(t) = A(t) + \rho^{-1}(1 - \rho) \left\{ A_G(t) - \Lambda_G^2(t) \right\}$. Clearly, $\Gamma_\rho(t) \leq \Gamma_1(t)$ if and only if $A_G(t) - \Lambda_G^2(t) \leq A(t) - \Lambda^2(t)$. Note that $A_G(t) - \Lambda_G^2(t) = \text{Var}\left\{ \int_0^t dN^c(s)/H(s) \right\}$ and $A(t) - \Lambda^2(t) = \text{Var}\left\{ \int_0^t dN^u(s)/H(s) \right\}$.

(4) Let $\rho \uparrow 1$, i.e., the proportion of missing $\delta_i$'s shrinks to 0. Then $\alpha^* \to 1$. The resulting estimator is $\hat{\Lambda}_1(t)$. On the other hand, if the censorship shrinks to 0, which entails $\Lambda_G(t) \to 0$ and $A_G(t) \to 0$, then $\alpha^* \to \rho$, which was the case discussed in the previous remark.

Table 3 displays the main results from our Monte Carlo studies on the adaptive estimator $\hat{F}(\hat{\alpha}^*, t) = 1 - e^{-\hat{\Lambda}(\hat{\alpha}^*, t)}$. The biases of the adaptive estimator and its variance estimator are small. The efficiency improvements of the adaptive estimator over the complete-case analysis and (to a lesser extent) over estimator (3.1) are impressive, especially for light censoring and substantial missingness.

*Proof of Theorem 3.1.* In analogy with the approximations given in Lemma 1, we can show that

$$(3.8) \quad \hat{\Lambda}_1(t) - \Lambda(t) = \frac{1}{n\rho} \left\{ \sum_{i=1}^n \int_0^t \frac{dM_i(s)}{H(s)}\xi_i + \sum_{i=1}^n \int_0^t \frac{Y_i(s)d\Lambda(s)}{H(s)}(\xi_i - \rho) \right.$$

$$\left. - \sum_{i=1}^n \Lambda(t)(\xi_i - \rho) \right\} + o_p(n^{-\frac{1}{2}})$$

$$= L_1(t) + o_p(n^{-\frac{1}{2}}), \text{ say,}$$

and

$$(3.9) \quad \hat{\Lambda}_2(t) - \Lambda(t) = \frac{1}{n(1-\rho)} \left[ \sum_{i=1}^n \int_0^t \frac{dM_i(s)}{H(s)}(1-\xi_i) - \sum_{i=1}^n \int_0^t \frac{Y_i(s)d\Lambda(s)}{H(s)}(\xi_i - \rho) \right.$$
$$\left. + \sum_{i=1}^n \Lambda(t)(\xi_i - \rho) - \sum_{i=1}^n \left\{ \int_0^t \frac{dN_i^c(s)}{H(s)} - \Lambda_G(t) \right\} \frac{\xi_i - \rho}{\rho} \right] + o_p(n^{-\frac{1}{2}})$$
$$= L_2(t) + o_p(n^{-\frac{1}{2}}), \text{ say.}$$

Thus to characterize the limiting distribution of $\hat{\Lambda}(\alpha, \cdot)$, it suffices to derive the covariance functions $E\{L_j(t)L_k(t')\}$ $(j, k = 1, 2)$. Through some tedious, but otherwise routine calculations, we obtain

$$(3.10) \qquad E\{L_1(t)L_1(t')\} = (n\rho)^{-1}\left\{A(t \wedge t') - (1-\rho)\Lambda(t)\Lambda(t')\right\},$$

$$(3.11) \qquad E\{L_1(t)L_2(t')\} = n^{-1}\left\{\Lambda(t)\Lambda(t') + \rho^{-1}\Lambda(t)\Lambda_G(t')\right\},$$

$$(3.12) \quad E\{L_2(t)L_2(t')\} = \{n(1-\rho)\}^{-1}\left[A(t \wedge t') + \rho^{-1}A_G(t \wedge t')\right.$$
$$\left. -\rho\left\{\Lambda(t) + \rho^{-1}\Lambda_G(t)\right\}\left\{\Lambda(t') + \rho^{-1}\Lambda_G(t')\right\}\right].$$

From (3.10)–(3.12), we can evaluate $nE\left[\{\alpha L_1(t) + (1-\alpha)L_2(t)\}\{\alpha L_1(t') + (1-\alpha)L_2(t')\}\right]$ to get the desired covariance function. □

We now return to the regression model studied in Section 2. Let $\hat{\beta}$ be as defined by (2.5). To estimate the cumulative baseline hazard function $\Lambda_0$, it is natural to extend the class of estimators given in (3.5). To avoid complicated asymptotic variances, we shall only consider $\alpha = 1$ and $\alpha = \hat{\rho}$, the two special cases discussed in Remarks (3) and (4) following Theorem 3.1. The two estimators for $\Lambda_0(t)$ are given below

$$(3.13) \qquad \hat{\Lambda}_1(\hat{\beta}, t) = \int_0^t \frac{\sum_{i=1}^n \xi_i dN_i^u(s)}{\hat{\rho}\sum_{i=1}^n Y_i(s)e^{\hat{\beta}'Z_i(s)}},$$

$$(3.14) \quad \hat{\Lambda}_2(\hat{\beta}, t) = \int_0^t \frac{\sum_{i=1}^n\{\xi_i dN_i^u(s) + (1-\xi_i)dN_i(s) - \hat{\rho}^{-1}(1-\hat{\rho})\xi_i dN_i^c(s)\}}{\sum_{i=1}^n Y_i(s)e^{\hat{\beta}'Z_i(s)}}.$$

**Theorem 3.2.** *Suppose that the assumptions of Theorem 2.3 are satisfied. Let $t_0 > 0$ be any number such that $EY_1(t_0) > 0$.*

(i) *The process $n^{1/2}\left\{\hat{\Lambda}_1(\hat{\beta}, \cdot) - \Lambda_0(\cdot)\right\}$ converges weakly in $\mathcal{D}[0, t_0]$ to a zero-mean Gaussian process with covariance function*

$$(3.15) \quad \tilde{\Gamma}_1(t, t') = \rho^{-1}\int_0^{t \wedge t'} \frac{d\Lambda_0(s)}{H_Z(s)} - \rho^{-1}(1-\rho)\Lambda_0(t)\Lambda_0(t')$$
$$+ a'(t)\Sigma(D)a(t') - \rho^{-1}(1-\rho)\left[a'(t)\Omega E\left\{N_1^{CZ}(\infty)\right\}\Lambda_0(t')\right.$$
$$\left. + a'(t')\Omega E\left\{N_1^{CZ}(\infty)\right\}\Lambda_0(t)\right],$$

where $H_Z(s) = E\left\{Y_1(s)e^{\beta_0'Z_1(s)}\right\}$, $a(t) = \int_0^t \tilde{z}(\beta_0, s)d\Lambda_0(s)$ *and* $\Omega = \{\rho V + (1 - \rho)DV\}^{-1}D$. *For fixed* $t$, $n^{1/2}\left\{\hat{\Lambda}_1(\hat{\beta}, t) - \Lambda_0(t)\right\} \xrightarrow{d} \mathcal{N}(0, \tilde{\Gamma}_1(t))$, *where*

$$(3.16) \qquad \tilde{\Gamma}_1(t) = \rho^{-1}\int_0^t \frac{d\Lambda_0(s)}{H_Z(s)} - \rho^{-1}(1 - \rho)\Lambda_0^2(t) + a'(t)\Sigma(D)a(t)$$
$$- 2\rho^{-1}(1 - \rho)a'(t)\Omega E\left\{N_1^{CZ}(\infty)\right\}\Lambda_0(t).$$

(ii) *The process* $n^{1/2}\{\hat{\Lambda}_2(\hat{\beta}, \cdot) - \Lambda_0(\cdot)\}$ *converges weakly to a zero-mean Gaussian process with covariance function*

$$(3.17)$$
$$\tilde{\Gamma}_2(t, t') = \int_0^{t \wedge t'} \frac{d\Lambda_0(s)}{H_Z(s)} + \rho^{-1}(1 - \rho)\text{Cov}\left\{N_1^{CH}(t), N_1^{CH}(t')\right\} + a'(t)\Sigma(D)a(t')$$
$$- \rho^{-1}(1 - \rho)\Bigg(a'(t)\Omega E\left[N_1^{CZ}(\infty)\left\{N_1^{CH}(t') - EN_1^{CH}(t')\right\}\right]$$
$$+ a'(t')\Omega E\left[N_1^{CZ}(\infty)\left\{N_1^{CH}(t) - EN_1^{CH}(t)\right\}\right]\Bigg),$$

*where* $N_i^{CH}(t) = \int_0^t dN_i^c(s)/H_Z(s)$ $(i = 1, \ldots, n)$. *For fixed* $t$, $n^{1/2}\{\hat{\Lambda}_2(\hat{\beta}, t) - \Lambda_0(t)\} \xrightarrow{d} \mathcal{N}(0, \tilde{\Gamma}_2(t))$, *where*

$$(3.18) \qquad \tilde{\Gamma}_2(t) = \int_0^t \frac{d\Lambda_0(s)}{H_Z(s)} + \rho^{-1}(1 - \rho)\text{Var}\left\{N_1^{CH}(t)\right\} + a'(t)\Sigma(D)a(t)$$
$$- 2\rho^{-1}(1 - \rho)a'(t)\Omega E\left[N_1^{CZ}(\infty)\left\{N_1^{CH}(t) - EN_1^{CH}(t)\right\}\right].$$

**Remarks.** (1) Consistent estimators for variances $\tilde{\Gamma}_1(t)$ and $\tilde{\Gamma}_2(t)$ may be obtained in a straightforward manner. For example, let $\hat{a}(t) = \int_0^t \bar{Z}(\hat{\beta}, s)d\hat{\Lambda}_1(\hat{\beta}, s)$ and $\hat{\Omega} = \left\{\hat{\rho}\hat{V} + (1 - \hat{\rho})\hat{V}D\right\}^{-1}D$. Then a consistent estimator for $\tilde{\Gamma}_1(t)$ is

$$n\hat{\rho}^{-1}\int_0^t \frac{d\hat{\Lambda}_1(\hat{\beta}, s)}{\sum_{i=1}^n Y_i(s)e^{\hat{\beta}'Z_i(s)}} - \hat{\rho}^{-1}(1 - \hat{\rho})\hat{\Lambda}_1^2(\hat{\beta}, t) + \hat{a}'(t)\hat{\Sigma}(D)\hat{a}(t)$$
$$- 2\hat{\rho}^{-1}(1 - \hat{\rho})\hat{a}'(t)\hat{\Omega}\left[\frac{1}{\sum_{i=1}^n \hat{\xi}_i}\sum_{i=1}^n \int_0^\infty \left\{Z_i(s) - \bar{Z}(\hat{\beta}, s)\right\}\xi_i dN_i^c(s)\right]\hat{\Lambda}_1(\hat{\beta}, t),$$

where $\hat{\Sigma}(D)$ is the consistent estimator given in Remark (2) following Theorem 2.3.

(2) If $D = 0$, then the last term on the right hand side of (3.16) disappears and the sum of the first and the third terms becomes the variance of the complete-case estimator. Thus, the use of $\hat{\Lambda}_1(\tilde{\beta}, t)$ reduces the variance by $\rho^{-1}(1 - \rho)\Lambda_0^2(t)$.

*Proof of Theorem 3.2.* Taking the Taylor expansions at $\beta_0$, we get, for $l = 1, 2$,

$$(3.19) \qquad \hat{\Lambda}_l(\hat{\beta}, t) = \hat{\Lambda}_l(\beta_0, t) - \int_0^t \bar{Z}'(\beta_0, s)d\hat{\Lambda}_l(\beta_0, s)(\hat{\beta} - \beta_0) + o_p(n^{-\frac{1}{2}}).$$

By the approximations given in the proofs of Theorems 2.1 and 2.2, we can express $\hat{\beta} - \beta_0$ approximately as a sum of $n$ i.i.d. random vectors. Furthermore, similar to (3.8),

$$
\begin{aligned}
(3.20) \qquad \hat{\Lambda}_1(\beta_0, t) - \Lambda_0(t) &= (n\rho)^{-1} \int_0^t \frac{\sum \xi_i dM_i(s)}{H_Z(s)} \\
&\quad + (n\rho)^{-1} \int_0^t \frac{\sum (\xi_i - \rho) \tilde{Y}_i(s)}{H_Z(s)} \lambda_0(s) ds \\
&\quad - (n\rho)^{-1} \Lambda_0(t) \sum_{i=1}^n (\xi_i - \rho) + o_p(n^{-\frac{1}{2}}) \\
&= J_1(t) + J_2(t) + J_3(t) + o_p(n^{-\frac{1}{2}}), \qquad \text{say},
\end{aligned}
$$

where $\tilde{Y}_i(s) = Y_i(s) e^{\beta_0' Z_i(s)}$. Let $\tilde{S}_k = \tilde{S}_{k1}(\infty) + \tilde{S}_{k2}(\infty)$ $(k = 1, 2)$, where $\tilde{S}_{kj}$ are defined in the proofs of Theorems 2.1 and 2.2. Then $E\left\{ \tilde{S}_1 J_3(t) \right\} = 0$ and

$$
\begin{aligned}
E\left\{ \tilde{S}_1 J_1(t) \right\} &= (1 - \rho) E \left\{ \int_0^\infty W_1(s) \tilde{Y}_1(s) d\Lambda_0(s) \int_0^{t \wedge s} \frac{dM_1(u)}{H_Z(u)} \right\} \\
&= -(1 - \rho) E \int_0^\infty \int_0^t 1_{(u \le s)} W_1(s) \tilde{Y}_1(s) H_Z^{-1}(u) \tilde{Y}_1(u) d\Lambda_0(s) d\Lambda_0(u).
\end{aligned}
$$

Moreover, we can show that $E\left\{ \tilde{S}_1 J_2(t) \right\} = -E\left\{ \tilde{S}_1 J_1(t) \right\}$. Therefore

$$
(3.21) \qquad E\left[ \tilde{S}_1 \left\{ J_1(t) + J_2(t) + J_3(t) \right\} \right] = 0.
$$

Likewise, we can show that $E\left[ \tilde{S}_{21}(\infty) \left\{ J_1(t) + J_2(t) + J_3(t) \right\} \right] = 0$. Thus

$$
\begin{aligned}
(3.22) \qquad & E\left[ \tilde{S}_2 \left\{ J_1(t) + J_2(t) + J_3(t) \right\} \right] \\
&= E\left[ \tilde{S}_{22}(\infty) \left\{ J_1(t) + J_2(t) + J_3(t) \right\} \right] \\
&= -\rho^{-1}(1 - \rho) E\left[ \left\{ N_1^{CZ}(\infty) - E N_1^{CZ}(\infty) \right\} \int_0^t \frac{dN_1^u(s)}{H_Z(s)} \right] \\
&= \rho^{-1}(1 - \rho) E N_1^{CZ}(\infty) \Lambda_0(t).
\end{aligned}
$$

From (3.21) and (3.22),

$$
\begin{aligned}
(3.23) \qquad & E\left[ \left\{ \tilde{S}_1 + D\hat{S}_2 \right\} \left\{ J_1(t) + J_2(t) + J_3(t) \right\} \right] \\
&= \rho^{-1}(1 - \rho) D E \left\{ N_1^{CZ}(\infty) \right\} \Lambda_0(t).
\end{aligned}
$$

It is also not difficult to show that

$$
\begin{aligned}
& E\left[ \left\{ J_1(t) + J_2(t) + J_3(t) \right\} \left\{ J_1(t') + J_2(t') + J_3(t') \right\} \right] \\
&= \frac{1}{n\rho} \int_0^{t \wedge t'} \frac{d\Lambda_0(s)}{H_Z(s)} - \frac{1 - \rho}{n\rho} \Lambda_0(t) \Lambda_0(t'),
\end{aligned}
$$

which, combined with (3.19), (3.20) and (3.23), yields the desired covariance function $\tilde{\Gamma}_1$.

For (ii), first note that

$$
\begin{aligned}
\hat{\Lambda}_2(\beta_0, t) &- \Lambda_0(t) \\
(3.24) \quad &= n^{-1} \sum_{i=1}^{n} \int_0^t \frac{dM_i(s)}{H_Z(s)} \\
&\quad - (n\rho)^{-1} \sum_{i=1}^{n} \left[ N_i^{CH}(t) - E\left\{ N_i^{CH}(t) \right\} \right] (\xi_i - \rho) + o_p(n^{-1/2}).
\end{aligned}
$$

From (3.24), we can show that $\left\{ \hat{\Lambda}_2(\beta_0, t) - \Lambda_0(t) \right\}$ is asymptotically uncorrelated with $\tilde{S}_1$ and $\tilde{S}_{21}$. The desired covariance formula (3.17) then follows by evaluating the asymptotic covariance between $\left\{ \hat{\Lambda}_2(\beta_0, t) - \Lambda_0(t) \right\}$ and $\tilde{S}_{22}$. The details are omitted. $\qquad\square$

## 4. Cox regression and cumulative hazard function estimation under Type II missingness

We now describe in detail the problem of Type II missingness mentioned in Section 1 using a slightly different notation. Let $(T_i^{(1)}, T_i^{(2)}, C_i, Z_i')$ $(i = 1, \ldots, n)$ be i.i.d. random vectors, where $T_i^{(1)}$ and $T_i^{(2)}$ denote two types of latent failure times, of which the first is of interest, and $C_i$ and $Z_i$ denote the censoring time and covariate vector as before. Suppose that, conditional on $Z_i$, the failure time $T_i^{(1)}$ is independent of $T_i^{(2)}$ and $C_i$, and has the hazard rate $\lambda(t \mid Z_i) = e^{\beta_0' Z_i} \lambda_0(t)$. Define $T_i = T_i^{(1)} \wedge T_i^{(2)}$, $\phi_i = 1_{(T_i^{(1)} \le T_i^{(2)})}$, $X_i = T_i \wedge C_i$ and $\delta_i = 1_{(T_i \le C_i)}$. Note that $\phi_i \delta_i$ indicates, by the value 1 vs. 0, whether or not the observation time $X_i$ is the failure time of interest $T_i^{(1)}$. In the standard competing risk setup, one observes $(X_i, \delta_i, \phi_i \delta_i, Z_i)$ for every $i$. With incomplete measurements on the failure types, however, the data consist of $(X_i, \delta_i, \xi_i, \xi_i \phi_i \delta_i, Z_i)$ $(i = 1, \ldots, n)$, where $\xi_i$ indicates, by the value 1 vs 0, whether $\phi_i$ is known or unknown. We assume that $\xi_i$ is independent of all other variables with $P(\xi_i = 1 \mid X_i, \delta_i, \phi_i, Z_i) = \tau$. This has the same level of generality as assuming $P(\xi_i = 1 \mid X_i, \delta_i = 1, \phi_i, Z_i) = \tau$, since for $\delta_i = 0$ the value $\xi_i$ does not have any effect on the observations and can therefore be redefined to make the independence true. We define $N_i^u$, $Y_i$, $\bar{Z}$ and $\bar{z}$ as in Section 2.

In the absence of missing values, the partial likelihood score function for $\beta_0$ is

$$
S^\phi(\beta) = \sum_{i=1}^{n} \int_0^\infty \left\{ Z_i(t) - \bar{Z}(\beta, t) \right\} \phi_i dN_i^u(t).
$$

By deleting all the cases with $\{ \delta_i = 1, \xi_i = 0 \}$, the complete-case estimating function is

$$
S_d^\phi(\beta) = \sum_{i=1}^{n} \int_0^\infty \left[ Z_i(t) - \frac{\sum_{j=1}^{n} \left\{ \delta_j \xi_j + (1 - \delta_j) \right\} Y_j(t) e^{\beta' Z_j(t)} Z_j(t)}{\sum_{j=1}^{n} \left\{ \delta_j \xi_j + (1 - \delta_j) \right\} Y_j(t) e^{\beta' Z_j(t)}} \right] \xi_i \phi_i dN_i^u(t).
$$

Because the index set $\{ j : (\delta_j \xi_j + (1 - \delta_j)) Y_j(t) = 1 \}$ is not a random subset of the risk set $\{ j : Y_j(t) = 1 \}$, the complete-case analysis does not yield a consistent estimator for $\beta_0$. We shall use the ideas presented in Section 2 to estimate $\beta_0$ under Type II missingness.

The analogs of estimating functions $S_1(\beta)$ and $S_2(\beta)$ studied in Section 2 are

$$S_1^\phi(\beta) = \sum_{i=1}^n \int_0^\infty \left\{ Z_i(t) - \bar{Z}(\beta, t) \right\} \xi_i \phi_i dN_i^u(t),$$

$$S_2^\phi(\beta) = \sum_{i=1}^n \int_0^\infty \left\{ Z_i(t) - \bar{Z}(\beta, t) \right\} \left\{ (1 - \xi_i) - \hat{\tau}^{-1}(1 - \hat{\tau})\xi_i(1 - \phi_i) \right\} dN_i^u(t),$$

where $\hat{\tau} = \sum_{i=1}^n \delta_i \xi_i / \sum_{i=1}^n \delta_i$. We have the following results for $S_k^\phi(\beta_0)$ $(k = 1, 2)$, which are similar to those of $S_k(\beta_0)$ $(k = 1, 2)$ given in Theorems 2.1 and 2.2.

**Theorem 4.1.** *The random vector* $n^{-1/2} \left[ S_1^\phi(\beta_0)', S_2^\phi(\beta_0)' \right]'$ *is asymptotically zero-mean normal with covariance matrix*

$$\begin{bmatrix} V_1^\phi & 0 \\ 0 & V_2^\phi \end{bmatrix},$$

*where* $V_1^\phi = \tau V^\phi$, $V_2^\phi = (1 - \tau)V^\phi + \tau^{-1}(1 - \tau)E(N_1^\phi - EN_1^\phi)^{\otimes 2}$, $N_i^\phi = \int_0^\infty \{ Z_i(t) - \bar{z}(\beta_0, t) \} (1 - \phi_i) dN_i^u(t)$ *and* $V^\phi = E \left[ \int_0^\infty \{ Z_1(t) - \bar{z}(\beta_0, t) \}^{\otimes 2} \phi_1 dN_1^u(t) \right]$.

In analogy with (2.5) for $\hat{\beta}$, we define $\hat{\beta}^\phi$ as a solution to

$$S_1^\phi(\beta) + D S_2^\phi(\beta) = 0,$$

where $D$ is a given $p \times p$ matrix. Then the following theorem similar to Theorem 2.3 holds.

**Theorem 4.2.** *Suppose that* $\left\{ \tau V^\phi + (1 - \tau)DV^\phi \right\}$ *is nonsingular. Then* $n^{1/2}(\hat{\beta}^\phi - \beta_0) \xrightarrow{d} \mathcal{N}(0, \Sigma^\phi(D))$, *where*

$$\Sigma^\phi(D) = \left\{ \tau V^\phi + (1 - \tau)DV^\phi \right\}^{-1} (\tau V^\phi + DV_2^\phi D') \left\{ \tau V^\phi + (1 - \tau)V^\phi D' \right\}^{-1}.$$

*The optimal choice for* $D$ *is* $D^* = (1 - \tau)V^\phi(V_2^\phi)^{-1}$, *in which case*

$$\Sigma^\phi(D^*) = \left\{ \tau V^\phi + (1 - \tau)^2 V^\phi(V_2^\phi)^{-1}V^\phi \right\}^{-1}.$$

*Proof of Theorem 4.1.* As in the proofs of Theorems 2.1 and 2.2, we can define the martingales $M_i^\phi(t) = \phi_i N_i^u(t) - \int_0^t Y_i(s)e^{\beta_0' Z_i(s)}\lambda_0(s)ds$ $(i = 1, \ldots, n)$ and derive the following key approximations

$$S_1^\phi(\beta_0) = \sum_{i=1}^n \int_0^\infty \left\{ Z_i(t) - \bar{z}(\beta_0, t) \right\} \xi_i dM_i^\phi(t)$$
$$+ \sum_{i=1}^n \int_0^\infty \left\{ Z_i(t) - \bar{z}(\beta_0, t) \right\} (\xi_i - \tau)Y_i(t)e^{\beta_0' Z_i(t)}\lambda_0(t)dt + o_p(n^{\frac{1}{2}}),$$

$$S_2^\phi(\beta_0) = \sum_{i=1}^n \int_0^\infty \left\{ Z_i(t) - \bar{z}(\beta_0, t) \right\} (1 - \xi_i)dM_i^\phi(t)$$
$$- \sum_{i=1}^n \int_0^\infty \left\{ Z_i(t) - \bar{z}(\beta_0, t) \right\} (\xi_i - \tau)Y_i(t)e^{\beta_0' Z_i(t)}\lambda_0(t)dt$$
$$- \tau^{-1} \sum_{i=1}^n (N_i^\phi - EN_i^\phi)(\xi_i - \tau) + o_p(n^{\frac{1}{2}}).$$

These two approximations can be used to show, through some tedious calculations, that the asymptotic variance-covariance matrix of $n^{-1/2} \left[ S_1^\phi(\beta_0)', S_2^\phi(\beta_0)' \right]'$ is $\text{Diag} \left\{ V_1^\phi, V_2^\phi \right\}$. Hence, the theorem follows from the multivariate central limit theorem. □

*Proof of Theorem 4.2.* This can be done by applying Theorem 4.1 and the arguments given in the proof of Theorem 2.3. □

Using $\hat{\beta}^\phi$ with its asymptotic distribution given by Theorem 4.2, we can construct consistent estimators for the cumulative baseline hazard function $\Lambda_0(t)$. Two such estimators which correspond to $\hat{\Lambda}_1(\hat{\beta}, t)$ and $\hat{\Lambda}_2(\hat{\beta}, t)$ defined by (3.13) and (3.14) are

$$\hat{\Lambda}_1^\phi(\hat{\beta}^\phi, t) = \int_0^t \frac{\sum_{i=1}^n \xi_i \phi_i dN_i^u(s)}{\hat{\tau} \sum_{i=1}^n Y_i(s) e^{\hat{\beta}^{\phi'} Z_i(s)}},$$

$$\hat{\Lambda}_2^\phi(\hat{\beta}^\phi, t) = \int_0^t \frac{\sum_{i=1}^n \left\{ \xi_i \phi_i + (1 - \xi_i) - \hat{\tau}^{-1}(1 - \hat{\tau}) \xi_i (1 - \phi_i) \right\} dN_i^u(s)}{\sum_{i=1}^n Y_i(s) e^{\hat{\beta}^{\phi'} Z_i(s)}}.$$

The kind of asymptotic properties given in Theorem 3.2 for $\hat{\Lambda}_k(\hat{\beta}, t)$ $(k = 1, 2)$ also hold for $\hat{\Lambda}_k^\phi(\hat{\beta}^\phi, t)$ $(k = 1, 2)$. They can be derived from the arguments used in the proof of Theorem 3.2. To simplify the statements, we shall only present the asymptotic normality part although the weak convergence also holds.

**Theorem 4.3.** *Suppose that the assumptions of Theorem 4.2 hold and that $t$ satisfies $EY_1(t) > 0$. Then $n^{1/2} \left\{ \hat{\Lambda}_k^\phi(\hat{\beta}^\phi, t) - \Lambda_0(t) \right\} \xrightarrow{d} \mathcal{N}(0, \sigma_{\phi,k}^2(t))$ $(k = 1, 2)$, where, with $H_Z(t)$ and $a(t)$ as defined in Theorem 3.2, $\Omega^\phi = \left\{ \tau V^\phi + (1 - \tau) DV^\phi \right\}^{-1} D$ and $N_i^{\phi H}(t) = \int_0^t (1 - \phi_i) dN_i^u(s)/H_Z(s)$ $(i = 1, \ldots, n)$,*

$$\sigma_{\phi,1}^2(t) = \tau^{-1} \int_0^t \frac{d\Lambda_0(s)}{H_Z(s)} - \tau^{-1}(1 - \tau)\Lambda_0^2(t) + a'(t)\Sigma^\phi(D)a(t)$$
$$- 2\tau^{-1}(1 - \tau)a'(t)\Omega^\phi EN_1^\phi \Lambda_0(t),$$

$$\sigma_{\phi,2}^2(t) = \int_0^t \frac{d\Lambda_0(s)}{H_Z(s)} + \tau^{-1}(1 - \tau)\text{Var} \left\{ N_1^{\phi H}(t) \right\} + a'(t)\Sigma(D)a(t)$$
$$- 2\tau^{-1}(1 - \tau)a'(t)\Omega^\phi E \left[ N_1^\phi \left\{ N_1^{\phi H}(t) - EN_1^{\phi H}(t) \right\} \right].$$

For the one-sample case, where the data consist of i.i.d. random vectors $(X_i, \delta_i, \xi_i, \xi_i \phi_i \delta_i)$ $(i = 1, \ldots, n)$, we modify (3.2), (3.4) and (3.5) to obtain the following class of consistent estimators for the cumulative hazard function of $T_i^{(1)}$:

$$\hat{\Lambda}^\phi(\alpha, t) = \alpha \int_0^t \frac{\sum_{i=1}^n \xi_i \phi_i dN_i^u(s)}{\hat{\tau} \sum_{i=1}^n Y_i(s)}$$
$$+ (1 - \alpha) \int_0^t \frac{\sum_{i=1}^n \left\{ 1 - \xi_i - \hat{\tau}^{-1}(1 - \hat{\tau}) \xi_i (1 - \phi_i) \right\} dN_i^u(s)}{(1 - \hat{\tau}) \sum_{i=1}^n Y_i(s)}.$$

The arguments given in the proof of Theorem 3.1 can be used to verify the following asymptotic normality for $\hat{\Lambda}^\phi(\alpha, t)$.

**Theorem 4.4.** *For $t$ satisfying $EY_1(t) > 0$, $n^{1/2} \left\{ \hat{\Lambda}^\phi(\alpha, t) - \Lambda(t) \right\} \xrightarrow{d} \mathcal{N}(0, \sigma_t^2(\alpha))$, where*

$$\sigma_t^2(\alpha) = \frac{\alpha^2}{\tau} \left\{ A(t) - (1-\tau)\Lambda^2(t) \right\} + 2\alpha(1-\alpha) \left\{ \Lambda^2(t) + \tau^{-1}\Lambda(t)\Lambda_Q(t) \right\}$$
$$+ \frac{(1-\alpha)^2}{1-\tau} \left[ A(t) + \tau^{-1} A_Q(t) - \tau \left\{ \Lambda(t) + \tau^{-1}\Lambda_Q(t) \right\}^2 \right],$$

*where $A(t) = \int_0^t d\Lambda(s)/EY_1(s)$, $A_Q(t) = \int_0^t d\Lambda_Q(s)/EY_1(s)$ and $\Lambda_Q$ is the cumulative hazard function of $T_i^{(2)}$. In particular,*

$$\sigma_t^2(1) = \tau^{-1} \left\{ A(t) - (1-\tau)\Lambda^2(t) \right\},$$
$$\sigma_t^2(\tau) = A(t) + \tau^{-1}(1-\tau) \left\{ A_Q(t) - \Lambda_Q^2(t) \right\}.$$

*The variance $\sigma_t^2(\alpha)$ is minimized when $\alpha$ equals*

$$\alpha^* = \frac{\tau \left\{ A(t) - \Lambda^2(t) \right\} + A_Q(t) - \Lambda_Q^2(t) - (1+\tau)\Lambda(t)\Lambda_Q(t)}{A(t) - \Lambda^2(t) + A_Q(t) - \Lambda_Q^2 - 2\Lambda(t)\Lambda_Q(t)}.$$

## 5. Discussions

We did not provide all the details for Type II missingness in Section 4 because of the similarity with Type I missingness. It should be noted that consistent estimators for the variance quantities such as $\Sigma^\phi(D)$, $\sigma_{\phi,k}^2(t)$ $(k = 1, 2)$ and $\sigma_t^2(\alpha)$ can be obtained in the same manners as their counterparts in Sections 2 and 3. Furthermore, the asymptotic approximations under Type II missingness are expected to have similar degrees of accuracy in finite samples as those of Type I missingness.

We have made the missing completely at random assumption in our developments. This assumption consists of two parts, the first part being the independence between $\xi_i$ and $(X_i, \delta_i, \phi_i, Z_i)$ for every $i$ and the second being the i.i.d property of $\xi_i$ $(i = 1, \ldots, n)$. The first part of the assumption cannot be avoided without direct modeling the missing processes. The second part can be relaxed to the extent that a consistent estimator for $P(\xi_i = 1)$ is available for every $i$. For example, in a multi-institutional study, it may be reasonable to assume that the missing probabilities are constant within the same institution but vary among different institutions. In this case, we may stratify our data on the institutions and modify the methods described in the previous sections to incorporate the stratification factor.

In many applications, the measurements on the covariate vectors are incomplete. [9] and subsequent papers provide solutions to this problem. It is possible to combine the techniques developed in this paper with those of [9] to handle the situation where both the covariates and the failure indicators are partially measured. The details will not be presented here.